\documentclass[a4paper,12pt,leqno,english]{article}
\usepackage[utf8]{inputenc}
\usepackage[T1]{fontenc}
\usepackage{babel}
\usepackage{amsmath}
\usepackage{amssymb}
\usepackage{hyperref}
\usepackage{tikz}

\newcommand{\C}{{\mathbb  C}}

\renewcommand{\L}{{\mathcal L}}

\newcommand{\N}{{\mathbb  N}}
\newcommand{\OO}{{\mathcal O}}

\newcommand{\Q}{{\mathbb  Q}}
\newcommand{\R}{{\mathbb  R}}

\newcommand{\scalar}[2]{{\langle#1,#2\rangle}}

\newcommand{\ch}{{\operatorname{ch}}}
\newcommand{\ext}{{\operatorname{ext}}}

\newcommand{\Log}{{\operatorname{Log}}}

\newcommand{\PSH}{{\operatorname{{\mathcal{PSH}}}}}

\newcommand{\supp}{{\operatorname{supp}\, }}

\newcommand{\LSC}{{\operatorname{{\mathcal{LSC}}}}}

\def\vfigura#1#2{
\setbox0\vbox{{
\input #1
}}
\setbox1\vbox{\hbox{\box0}\hbox{{\obeylines #2}}}
\dimen0 = -\ht1
\advance\dimen0 by-\dp1
\dimen1 = \wd1
\dimen2 = -\dimen0
\divide\dimen2 by\baselineskip
\count100 = 1
\advance\count100 by\dimen2
\advance\count100 by1
\box1
\hangindent\dimen1
\hangafter=-\count100
\vskip\dimen0
}

\numberwithin{equation}{section}

\newtheorem{theorem+}           {Theorem}      [section]
\newtheorem{definition+}  [theorem+]  {Definition}
\newtheorem{lemma+}  [theorem+]  {Lemma}
\newtheorem{corollary+}  [theorem+]  {Corollary}
\newtheorem{proposition+}  [theorem+]  {Proposition}
\newtheorem{example+}  [theorem+]  {Example}
\newtheorem{problem+}  [theorem+]  {Problem}
\newtheorem{remark+}  [theorem+]  {Remark}

\newenvironment{theorem}{\begin{theorem+}\sl}{\end{theorem+}\rm}

\newenvironment{corollary}{\begin{corollary+}\sl}{\end{corollary+}\rm}
\newenvironment{proposition}{\begin{proposition+}\sl}{\end{proposition+}\rm}

\newenvironment{proof}{\medbreak\noindent{\bf  Proof:}\rm}{\hfill$\square$\rm}
\newenvironment{prooftx}[1]{\medbreak\noindent{\bf 
#1:}\rm}{\hfill$\square$\rm} 

\title{{\Large \bf 
Polynomials with exponents in  compact convex sets  and 
associated weighted extremal functions -  \\
The Bernstein-Walsh-Siciak theorem}}
\author{{Benedikt Steinar Magnússon, Ragnar Sigurðsson} \\
{and Bergur Snorrason}}
\date{{}}

\begin{document}

\maketitle

\begin{abstract} \noindent
We generalize the Bernstein-Walsh-Siciak theorem  on polynomial approximation
in $\C^n$ to the case  where the polynomial ring 
${\mathcal  P}(\C^n)$ is replaced by a subring ${\mathcal P}^S(\C^n)$ 
consisting of all polynomials with exponents restricted to 
sets $mS$, where $S$ is a compact convex subset of $\R^n_+$
with $0\in S$ and $m=0,1,2,3,\dots$, and uniform estimates of
error in the approximation are replaced by weighted uniform estimates 
with  respect to an admissible weight function.

\end{abstract}


\section{Introduction}
\label{sec:01}

\noindent
The  Runge-Oka-Weil theorem states that if $K$ is a compact
polynomially convex subset of $\C^n$ and $f$ is a holomorphic function
in some neighborhood of $K$, then $f$ can be approximated uniformly on
$K$ by polynomials.  We let ${\mathcal P}_m(\C^n)$ denote the space
of polynomials of degree $\leq m$ in $n$ complex variables and let 
$$
d_{K,m}(f)=\inf\{\|f-p\|_K\,;\, p\in {\mathcal P}_m(\C^n)\}
$$
denote the smallest error in an approximation of $f$ by polynomials
of degree $\leq m$, i.e., the distance from $f$ to 
${\mathcal P}_m(\C^n)$ in the supremum norm $\|\cdot\|_K$ on 
${\mathcal C(}K)$.  Then the Runge-Oka-Weil theorem 
is equivalent to stating that
$$
\lim_{m\to \infty}d_{K,m}(f)=0. 
$$
The Bernstein-Walsh-Siciak theorem states that 
$$
\varlimsup_{m\to \infty} d_{K,m}(f)^{1/m}\leq \dfrac 1R
$$
if and only if $f$ has a holomorphic extension to 
$X_R=\{z\in \C^n\,;\, \Phi_K(z)<R\}$,  
where $\Phi_K=\varlimsup_{m\to\infty}\Phi_{K,m}$ and 
$\Phi_{K,m}=\sup\{|p|^{1/m}\,;\, p\in {\mathcal P}_m(\C^n), \|p\|_K\leq 1\}$
are the Siciak functions of the set $K$ and it is  assumed that
$\Phi_K$ is continuous and $R\geq 1$.  
For more on this result see Siciak \cite[\S 10]{Sic:1962}.

\smallskip
In this paper, we  generalize  the Bernstein-Walsh-Siciak theorem
where   ${\mathcal P}_m(\C^n)$ is replaced by
the space ${\mathcal P}^S_m(\C^n)$ of  all polynomials $p$ of the form 
$$
p(z)=\sum_{\alpha\in (mS)\cap \N^n} a_\alpha z^\alpha, \quad z\in \C^n,
$$
for a given compact convex subset $S$ of $\R^n_+$ with $0\in S$,
and ${\mathcal P}(\C^n)$ is replaced by
${\mathcal P}^S(\C^n)=\cup_{m\in \N}{\mathcal P}^S_m(\C^n)$.
For every function $q\colon E\to \R\cup\{+\infty\}$
defined on a subset $E$ of $\C^n$  and
$m=1,2,3,\dots$ we define 
the \emph{Siciak functions with respect to $S$, $E$, $q$, and $m$} by
$$
\Phi^S_{E,q,m}(z)=\sup\{|p(z)|^{1/m}\,;\, p\in {\mathcal P}^S_m(\C^n),
\|pe^{-mq}\|_E\leq 1\},
\quad
z \in \mathbb{C}^n,
$$   
and the Siciak function with respect to $S$, $E$, and $q$   by
$\Phi^S_{E,q}=\varlimsup_{m\to\infty}\Phi^S_{E,q,m}$. 
By Proposition 2.2 in \cite{MagSigSigSno:2023} we have that 
\begin{equation}
\label{eq:1.1}
    \Phi^S_{E,q}(z)
    =
    \sup_{m \in \mathbb{N}} \Phi^S_{E,q,m}(z)
    =
    \lim_{m\to\infty} \Phi^S_{E,q,m}(z),
    \quad
    z \in \mathbb{C}^n,
\end{equation}
and if $q$ is bounded below and $\Phi^S_{E, q}$ is continuous on a compact $X \subset \mathbb{C}^n$,
    then the convergence is uniform on $X$.

\smallskip
The Lelong class with respect to $S$, denoted by $\L^S(\C^n)$, is
defined in terms of the supporting function $\varphi_S$ of $S$, given by
$\varphi_S(\xi)=\sup_{s\in S}\scalar s\xi$ for  $\xi\in \R^n$,
and the map  $\Log\colon \C^{*n}\to \R^n$,
given by $\Log\,  z=(\log|z_1|,\dots,\log|z_n|)$.    
We define the \emph{logarithmic supporting function}
$H_S\in \PSH(\C^n)$ of $S$ by 
$$H_S=\varphi_S\circ \Log$$
on $\C^{*n}$ and extend the definition to the coordinate 
hyperplanes  by the formula
$$
H_S(z)=\varlimsup_{\C^{*n}\ni w\to z}H_S(w), \quad z\in \C^n\setminus
\C^{*n}.
$$
The Lelong class $\L^S(\C^n)$ 
is defined as the set of all $u\in \PSH(\C^n)$ satisfying a growth
estimate of the form $u\leq c_u+H_S$ for some constant $c_u$. 
The \emph{Siciak-Zakharyuta function with respect to $S$, $E \subset \mathbb{C}^n$, and
    $q \colon E \to \R \cup \{+\infty\}$} is defined by 
$$
V^S_{E,q}(z)=\sup\{u(z)\,;\, u\in \L^S(\C^n), u|_E\leq q\},
\quad
z \in \mathbb{C}^n.
$$ 
Finally, the distance from a bounded function $f$ on $E$   
to ${\mathcal P}^S_m(\C^n)$
with respect to the supremum norm on $E$ with weight $e^{-mq}$
is defined by
$$
d^S_{E,q,m}(f)=\inf\{\|(f-p)e^{-mq}\|_E\,;\, p\in {\mathcal P}^S_m(\C^n)\}.
$$ 
We call the sequence $(d_{E,q,m}^{S}(f))_{m\in \N}$ the \emph{approximation numbers
    of $f$ on $E$ with respect to $S$ and $q$}.

Observe that the standard simplex
$\Sigma=\ch\{0,e_1,\dots,e_n\}$
has the supporting function
$\varphi_\Sigma(\xi)=\max\{\xi_1^+\dots,\xi_n^+\}$
where $\xi_j^+=\max\{\xi_j,0\}$ and that the 
logarithmic supporting function is 
$H_\Sigma(z)=\log^+\|z\|_\infty$, which implies that
the Lelong class
$\L^\Sigma(\C^n)$ is the standard Lelong class $\L(\C^n)$, see for example Section $5$ in Klimek \cite{Kli:1991},
and the polynomial space 
${\mathcal P}^\Sigma_m(\C^n)$ is ${\mathcal P}_m(\C^n)$.
We drop the superscript $S$ in the case $S=\Sigma$ and
the subscript $q$ in the case $q=0$.  

\smallskip
The function $q \colon E \rightarrow \mathbb{R} \cup \{+\infty\}$ is said to be an
\emph{admissible weight  with respect to $S$
on $E$} if  $q$ is lower semi-continuous,   
the set $\{z\in E \,;\, q(z)<+\infty\}$ is non-pluripolar,
  and  for  $E$ unbounded
$\lim_{E\ni z, |z|\to \infty}(H_S(z)-q(z))=-\infty$.

\smallskip
By Proposition 3.6 in \cite{MagSigSigSno:2023},
an entire function $p$ is in  ${\mathcal
  P}^S_m(\C^n)$ if and only if $\log|p|^{1/m}$ is in $\L^S(\C^n)$. 
This implies that $\log \Phi^S_{E,q}\leq V^S_{E, q}$. 
In Corollary 4.7 in \cite{MagSigSigSno:2023},
we have  examples where equality does not hold and in 
Theorem 1.1 in \cite{MagSigSig:2023} it is proved that for an
admissible weight $q$ on a closed set $E \subset \mathbb{C}^n$ the equality
$V^S_{E,q}=\log\Phi^S_{E,q}$ holds on $\mathbb{C}^{*n}$ if and only if $S\cap \Q^n$ is 
dense in  $S$.

\smallskip
Bos and Levenberg \cite[Theorem 3.1]{BosLev:2018} generalized 
the Bernstein-Walsh-Siciak theorem to the case
where the polynomial ring is ${\mathcal P}^S(\C^n)$ 
and the weight is $q=0$, with the assumption that
$V^S_K$ is continuous and $S$ is a lower set.
We say that $S$ is a \emph{lower set} if for all $x \in S$
    the box $[0, x_1] \times \dots \times [0, x_n] \subset S$.
Lower sets can be described in terms of their supporting functions,
    logarithmic supporting functions, and
    respective Lelong classes.
See Theorem 5.8 in \cite{MagSigSigSno:2023} and Theorem 5.1 in \cite{Sno:2024b}.

\smallskip 
We have not been able to find a generalization of the Bernstein-Walsh-Siciak theorem,
with weighted uniform estimates, in the form of an equivalence 
statement.  We have to separate our results into two
parts, where in the first part we assume that the 
approximation numbers $d^S_{K,q,m}(f)$ decrease 
exponentially with $m$,
    and in the second we assume $f$ can be extended to a holomorphic function
    on a sublevel set of $V^S_{K, q}$.
One reason for this separation is the different influences $S$ and $q$ have.
The added detail to Theorem \ref{thm:1.1} is due to the weight $q$,
    while most of the added assumptions in Theorem \ref{thm:1.2} are because of $S$.
Point {\bf (iii)} in Theorem \ref{thm:1.1} and Corollary \ref{cor:1.3} give sufficient
    conditions for these added assumptions to be superfluous.

\begin{theorem}
  \label{thm:1.1}
Let $S \subset \R_+^n$ be a compact convex set
with $0\in S$, $q$ be an admissible weight on
a compact subset $K$ of $\C^n$,
such that $V^{S*}_{K,q} \leq q$ on $K$, and  
for every $r>0$ define  
  $X_{r}=\{ z\in \C^n\,;\, 
V^S_{K,q}(z) < \log r\}$.
Let  $f\colon K\to \C$ be bounded, assume that
\begin{equation*}
L=\{z\in K\,;\, \lim_{m\to \infty} d_{K,q,m}^{S}(f)e^{mq(z)}=0 \}
\neq \varnothing
\end{equation*}
and that  
for some $R>0$,  $K\subset X_{R}$,  and
\begin{equation}
  \label{eq:1.2}
  \varlimsup_{m\to \infty} d^S_{K,q,m}(f)^{1/m}\leq \dfrac 1R.
\end{equation}
Then the following hold:
\begin{itemize}
    \item[{\bf (i)}]
        For every $\gamma>0$, such that $K\subset X_{R-\gamma}$, the function  $f|_L$ extends to a holomorphic function
            $F_\gamma\in \OO(X_{R-\gamma})$.
    \item[{\bf (ii)}]
        If $X$ is an open component of $X_{R}$ and
            $L_X=L\cap X$ is non-pluripolar,
            then $f|_{L_X}$ extends to a unique holomorphic function on $X$. 
    \item[{\bf (iii)}]
        If $q < \log R$ then $L = K$ and, consequently,
            there exists $F \in \mathcal{O}(X_R)$ such that $F|_K = f$.
\end{itemize}
\end{theorem}

If $p_m$ are such that $\|(f - p_m)e^{-mq}\|_K$ is close to $d^S_{K, q, m}(f)$,
    as in Section \ref{sec:02}, then $L$ denotes the set where $p_m \to f$, pointwise.
If $q$ is bounded above then $p_m \to f$ uniformly on $L$.

The assumption that $V^{S*}_{K, q} \leq q$ on $K$ implies that $V^S_{K, q}$ is upper semicontinuous
    and that the sublevel sets $X_r$ are all open.
See Proposition $4.5$ in \cite{MagSigSigSno:2023}.
Here $V^{S*}_{K, q}(z) = \varlimsup_{w \rightarrow z} V^S_{K, q}(w)$
    denotes the \emph{upper semicontinuous regularization} of $V^S_{K, q}$.

\medskip
For the converse statement, we need the concept of a hull 
of a compact convex subset $S$ of $\R^n_+$ with respect
to a cone $\Gamma$,
\begin{equation}
  \label{eq:1.3}
\widehat S_\Gamma=\{x\in \R^n_+ \,;\,
\scalar x\xi \leq \varphi_S(\xi), \forall \xi\in \Gamma\}.  
\end{equation}
We have  
$S=\widehat S_{\R^n}$ for every compact convex subset $S$ of $\R^n_+$
and if $\Gamma_1\subseteq \Gamma_2$ then 
$\widehat S_{\Gamma_2}\subseteq \widehat S_{\Gamma_1}$.
Since $S \subset \mathbb{R}^n_+$ and $0\in S$, we have that $\varphi_S=0$ on $\R^n_-$ and 
$S=\widehat S_{(\R^n\setminus \R^n_-)\cup\{0\}}$.

\begin{theorem}  \label{thm:1.2}
Let $S\subset \R_+^n$ be a compact convex set with $0\in S$ and $S = \overline{S \cap \mathbb{Q}^n}$,
    $q$ be an admissible weight  on 
a compact subset $K$  of $\C^n$, and $R>0$.
Assume that $V^{S}_{K,q}$ is continuous, $X_R$ is bounded, and 
$K\subset X_R$.  
Let $d_m=d(mS,\N^n\setminus mS)$ be  the euclidean distance between
$mS$ and $\N^n\setminus mS$.
If $f\in \OO(X_R)$ then
\begin{equation}
  \label{eq:1.4}
  \varlimsup_{m\to \infty} d^{\widehat S_{\Gamma_m}}_{K,q,m}(f)^{1/m}\leq 
\dfrac 1{R},
\end{equation}
where 
$\Gamma_m
=\{\xi\in \R^n\,;\,  \scalar{\mathbf 1}\xi\geq -\tfrac 12 d_m|\xi| \},   
$
and the hull $\widehat S_{\Gamma_m}$ with respect to $\Gamma_m$
is defined by $(\ref{eq:1.3})$.
\end{theorem}

In the special case when $S$ is a lower set
    we have $\widehat S_{\R^n_+}=S$ by Theorem 5.8 in \cite{MagSigSigSno:2023}.
So $S \subset \widehat{S}_{\Gamma_m} \subset \widehat{S}_{\mathbb{R}^n_+} = S$,
    since $\mathbb{R}^n_+ \subset \Gamma_m$.
It is possible that $S$ is a lower set, but not a convex body,
    that is $S$ may be a lower set with an empty interior.
These lower sets, however, are not interesting in their own right,
    since they can be written, after possibly rearranging the variables, as
    $S = T \times \{0\}^{n - \ell}$, where $T \subset \mathbb{R}^\ell_+$ is lower body.
In this case $\varphi_S(\xi) = \varphi_T(\xi_1, \dots, \xi_\ell)$, so the function
    in $\mathcal{L}^S(\mathbb{C}^n)$ only depend on the first $\ell$ variables.
When $S$ is a lower body we also have that the sublevel sets $X_r$ are bounded, for $r > 0$.
This holds since there exists $a > 0$ such that $a\Sigma \subset S$, so $H_S \geq a \log^+\| \cdot \|_\infty$.
By Proposition 4.5 in \cite{MagSigSigSno:2023}, there exists a constant $c$ such that
    $V^S_{K, q} \geq H_S + c \geq a\log\| \cdot \|_\infty + c$.
So when $S$ is a lower body the sublevel sets of $V^S_{K, q}$ are bounded.

\begin{corollary}
  \label{cor:1.3}
Let $S$ be a lower set,
 $q$ be an admissible weight on a compact subset $K$ of $\C^n$ and 
assume that $V^S_{K,q}$ is continuous and $K\subset X_R$ for some
$R>0$.
If $f\in \OO(X_R)$, then \eqref{eq:1.2} holds.
\end{corollary}

A Bernstein-Walsh-Siciak theorem for a lower set $S$ and
    weight $q = 0$ is proved in Bos and Levenberg \cite[Theorem 3.1]{BosLev:2018}. 
Their result follows from Corollary \ref{cor:1.3}, since if $q = 0$ then $K \subset X_R$ if and only if $R > 1$.

\begin{corollary}
  \label{cor:1.4}
Let $S$ be a lower set, $K$ be a compact subset of 
$\C^n$, $R>1$, and assume that $V^S_{K}$ is continuous.
Then  $f\in \OO(X_R)$ if and only if $\varlimsup\limits_{m \to \infty} d^S_{K, m}(f)^{1/m} \leq 1/R$.
\end{corollary}

Setting $S = \Sigma$ in Corollary 1.3 we get a weighted Bernstein-Walsh-Siciak theorem.
We include it to justify our definition of the approximation numbers $d^S_{K, q, m}(f)$.
In \cite[Theorem 2.4]{ChaLevWie:2024} a weighted Bernstein-Walsh theorem is proved, for $n = 1$,
    using approximation by weighted polynomials.
They use $d'_{K, q, m}(f) = \inf\{\|f - e^{-mq}p\| \,;\, p \in \mathcal{P}_m(\mathbb{C}^n)\}$,
    with the assumption that $e^{-q}$ extends as an entire function which is non-vanishing on $K$.

Approximations of this form have been studied for decades,
    but it does not match some of our intuitions for $V^S_{K, q}$.
To see this we note that $d'_{K, q, m}(f) = d_{K, m}(f)$ when $q = a$ is any constant weight,
    whereas $d^S_{K, q, m}(f) = e^{-ma}d^S_{K, m}(f)$ under the same assumption.
In this case, $d^S_{K, q, m}(f)$ are a closer match for the classical Bernstein-Walsh-Siciak theorem,
    since $V^S_{K, q} = V^S_K + a$.

\begin{corollary}
\label{cor:1.5}
Let $R > 0$ and $q$ be an admissible weight on a compact subset $K$ of $\C^n$, and 
assume that $V^*_{K,q} \leq q$ and $K\subset X_R$.
If $f\in \OO(X_R)$, then $\varlimsup\limits_{m \to \infty} d_{K, q, m}(f)^{1/m} \leq 1/R$.
\end{corollary}

Section \ref{sec:04} discusses lower bounds for the distances from $mS$ to $\mathbb{N}^n \setminus mS$.
This is done to justify some technical decisions in the proof of Theorem \ref{thm:1.2}.
They are also interesting in their own right
    since these distances appear in results that give sufficient conditions for entire functions to belong
    to certain polynomial classes.
See Theorems $3.6$ and $7.2$ in \cite{MagSigSigSno:2023}.

This paper is part of a series of papers studying the Lelong class $\mathcal{L}^S(\mathbb{C}^n)$,
    with the aim of relating it to polynomial approximations using polynomials from $\mathcal{P}^S(\mathbb{C}^n)$.
This self-contained exposition began in \cite{MagSigSigSno:2023},
    with a focus on fundamental results that could be useful in the following papers.

A commonly used tool in pluripotential theory is approximation by integral convolution with a smoothing kernel.
The Lelong class $\mathcal{L}(\mathbb{C}^n)$ is closed under such smoothing,
    whereas $\mathcal{L}^S(\mathbb{C}^n)$ may not be.
In \cite[Theorem 5.8]{MagSigSigSno:2023} it is shown that $\mathcal{L}^S(\mathbb{C}^n)$ is closed under such smoothing
    if and only if $S$ is a lower set.
Methods of approximation under which $\mathcal{L}^S(\mathbb{C}^n)$ is always closed are considered in \cite{Sno:2024a}.

Siciak \cite[Proposition 5.9]{Sic:1981} proved a product formula for the Siciak-Zakharyuta functions,
    relating the behavior
    of the Siciak-Zakharyuta function of a cartesian product with the Siciak-Zakharyuta function of each term.
This formula is greatly generalized in \cite{Sno:2024b}, along with showing that a weighted version of the formula
    is not possible.

The Siciak-Zakharyuta theorem relates the Siciak-Zakharyuta functions with the Siciak functions.
It states that $V_K = \log \Phi_K$.
In \cite{MagSigSig:2023} it is showed that $V^S_{K, q} = \log \Phi^S_{K, q}$ on $\mathbb{C}^{*n}$ if and only if
    $S \cap \mathbb{Q}^n$ is dense in $S$.
This result will play an important role in the proof of Theorem \ref{thm:1.2}.

In \cite{Sig:2024} generalizations of the Runge-Oka-Weil theorem are considered.
Namely, the paper considers when holomorphic functions can be approximated by polynomials in
    $\mathcal{P}^S(\mathbb{C}^n)$.
The hull of compact subsets of $\mathbb{C}^n$ with respect to $\mathcal{P}^S(\mathbb{C}^n)$ is relevant,
    as in the classical Runge-Oka-Weil theorem.
This hull is the same as the classical polynomial hull when $S$ contains a neighborhood of $0$ in $\mathbb{R}^n_+$,
    but is otherwise more complicated, as it may not be bounded.
Techniques from \cite{Sig:2024} appear in the proof of Theorem \ref{thm:1.2}.

\subsection*{Acknowledgment}
The results of this paper are a part of a research project, 
\emph{Holomorphic Approximations and Pluripotential Theory},
with  project grant 
no.~207236-051 supported by the Icelandic Research Fund.
We would like to thank the Fund for its support and the 
Mathematics Division, Science Institute, University of Iceland,
for hosting the project.   We thank Álfheiður Edda Sigurðardóttir
for many fruitful discussions and a careful reading of the manuscript.

\section{Proof of Theorem \ref{thm:1.1}}
\label{sec:02}
The proof is a modification of a classical
argument which is the easy part of the equivalence in the original 
Bernstein-Walsh theorem. 

\begin{prooftx}{Proof of Theorem \ref{thm:1.1}}
For a sequence $\varepsilon_m \searrow 0$ we can, by \eqref{eq:1.2},
    find a sequence of polynomials $(p_m)_{m \in \mathbb{N}}$,
    with $p_0 = 0$ and $p_m \in \mathcal{P}^S_m(\mathbb{C}^n)$,
    such that
\begin{equation}
\label{eq:2.1}
  |f(z)-p_m(z)|\leq (1+\varepsilon _m)
d_{K,q,m}^{S}(f) e^{mq(z)}, \quad z\in K.
\end{equation}
Then, for every $\gamma \in ]0, R[$ such that $K \subset X_{R - \gamma}$,
there exists a constant $A_\gamma>0$, such that
\begin{equation*}
d^S_{K,q,m}(f)\leq  \|(f-p_m)e^{-mq}\|_K\leq
\dfrac{A_\gamma(1+\varepsilon_m)^m}{(R-\gamma)^m}, \quad m\in \N.
\end{equation*}
For every $j\in\N^*$ and every $z\in K$ we have
\begin{align*}
  |p_j(z)-p_{j-1}(z)| &\leq |f(z)-p_j(z)|+ |f(z)-p_{j-1}(z)| \\
&\leq 
\dfrac{A_\gamma (1+\varepsilon_{j-1})^{j}}{(R-\gamma)^j} 
\Big(1+\dfrac {R-\gamma}{e^{q(z)}}\Big)\cdot e^{jq(z)}. 
\end{align*}
Since $q\in \LSC(K)$ takes its minimum $a$ at some point in 
$K$, we have 
\begin{equation*}
  \dfrac 1j
  \log\big((R-\gamma)^j|p_j(z)-p_{j-1}(z)|/(1+\varepsilon_{j-1})^{j}
B_\gamma\big)
\leq q(z), \quad z\in K,
\end{equation*}
where $B_\gamma=A_\gamma(1+(R-\gamma)/e^{a})$.  By the
definition of $V^S_{K,q}$ this implies that
\begin{equation*}
  \dfrac 1j
  \log\big((R-\gamma)^j|p_j(z)-p_{j-1}(z)|/(1+\varepsilon_{j-1})^{j}
B_\gamma\big)
\leq V^S_{K,q}(z), \quad z\in \C^n,
\end{equation*}
and consequently 
\begin{equation*}
  \big|p_j(z)-p_{j-1}(z)\big| \leq \dfrac{B_\gamma
    (1+\varepsilon_{j-1})^{j}
e^{jV^S_{K,q}(z)}}
{(R-\gamma)^j}, \quad z\in \C^n.
\end{equation*}
If $0 < \varrho < 1$ such that $K \subset X_{\varrho(R - \gamma)}$,
then
$V^S_{K,q}(z)\leq \log(\varrho(R-\gamma))$ for $z\in
X_{\varrho(R-\gamma))}$, which implies that  
\begin{equation*}
  |p_j(z)-p_{j-1}(z)|\leq B_\gamma
  \big((1+\varepsilon_{j-1})\varrho\big)^j, 
\quad z\in 
X_{\varrho(R-\gamma)},
\end{equation*}
and this estimate implies that $p_m=\sum_{j=1}^m(p_j-p_{j-1})$
converges locally uniformly on $X_{R-\gamma}$ to a holomorphic
function $F_\gamma$.
By \eqref{eq:2.1},
    we have  $F_\gamma=f$ on $L$, proving {\bf (i)}.
Point {\bf (ii)} then follows from the identity principle for holomorphic functions.  
For {\bf (iii)} we assume $\sup_K q < \log R$ and let $z \in K$ and $R_0 > 0$ such that $q < \log R_0 < \log R$.
Then $e^{q(z)} < R_0$, so
\begin{equation*}
    \varlimsup_{m \rightarrow \infty}
    \left (
        d^S_{K, q, m}(f)e^{mq(z)}
    \right )^{1/m}
    <
    R_0
    \varlimsup_{m \rightarrow \infty}
        d^S_{K, q, m}(f)^{1/m}
    \leq
    \frac{R_0}{R}
    <
    1.
\end{equation*}
Therefore, $\lim_{m \rightarrow \infty} d^S_{K, q, m}(f)e^{mq(z)} = 0$,
    since $\sum_{m = 0}^\infty d^S_{K, q, m}(f)e^{mq(z)}$ converges by the root test.
Since this holds for all $z \in K$ we have that $L = K$.
\end{prooftx}

\section{Proof of Theorem \ref{thm:1.2}}
Proving \eqref{eq:1.4} involves constructing polynomials that approximate $f$ sufficiently well.
This is done by solving $\bar \partial$-equations using Hörmander's $L^2$-estimates.

\begin{prooftx}{Proof of Theorem \ref{thm:1.2}}
It is sufficient to 
construct a sequence of polynomials $(p_m)_{m \in \mathbb{N}}$ such that
    $p_m \in \mathcal{P}^{\widehat S_{\Gamma_m}}_m(\mathbb{C}^n)$ and 
    \begin{equation}
      \label{eq:3.1}
      \varlimsup_{m\to \infty}\|(f - p_m)e^{-mq}\|_K^{1/m}\leq 1/R.
    \end{equation}
The polynomial $p_m$ will be of the form $p_m = \chi f - u_m$ where
$\chi$ is a cut-off function with 
    support in $X_R$,
    equal to $1$ on $\overline{X}_{R - \gamma}$, for some $\gamma \in
    ]0, R[$, and 
    $u_m$ is a solution to $\bar \partial u_m = f \bar \partial \chi$
    satisfying the weighted $L^2$ estimate 
\begin{equation}
\label{eq:3.2}
    \int_{\mathbb{C}^n}
        |u_m|^2 (1 + |w|^2)^{-a_m} e^{-2mV_m}\, d\lambda
    \leq
    \frac{1}{a_m}
    \int_{\mathbb{C}^n}
        |f \bar \partial \chi|^2 (1 + |w|^2)^{-a_m + 2} e^{-2mV_m}\, d\lambda,
\end{equation}
where $(V_m)_{m \in \mathbb{N}}$ is a sequence in $\mathcal{L}^S(\mathbb{C}^n)$ and
    $(a_m)_{m \in \mathbb{N}}$ is a sequence of strictly positive numbers such that
\begin{equation}
\label{eq:3.3}
    \varliminf_{m \rightarrow \infty} a_m^{1/m}
    =
    1.
\end{equation}
For the choice of these sequences we follow Section $5$ in Sigurðardóttir \cite{Sig:2024}.

Since $V^S_{K,q}$ is continuous, 
the Siciak-Zakharyuta theorem, Theorem $1.1$ in \cite{MagSigSig:2023},
implies that $\log \Phi^S_{K,q} = V^S_{K,q}$ and by \eqref{eq:1.1} we have that
$\log \Phi^S_{K,q,m} \rightarrow V^S_{K,q}$ 
locally uniformly on $\mathbb{C}^n$.

The subset $S_m = \ch(S \cap \frac{1}{m} \mathbb{N}^n)$ of $S$ is a
polytope with rational vertices and 
    $\mathcal{P}^{S_m}_m(\mathbb{C}^n) = \mathcal{P}^S_m(\mathbb{C}^n)$.
Hence $\Phi^{S_m}_{K,q,m} = \Phi^S_{K,q,m}$, and
\begin{equation*}
    \log \Phi^S_{K,q,m}
    =
    \log \Phi^{S_m}_{K,q,m}
    \leq
    V^{S_m}_{K,q}
    \leq
    V^S_{K,q}.
\end{equation*}
Consequently, $V^{S_m}_{K,q} \rightarrow V^S_{K,q}$ uniformly on compact
subsets of $\mathbb{C}^n$. 
For simplicity, we set 
$a_m = \frac{1}{2} d(mS_m, \mathbb{N}^n
\setminus mS_m)$, $V_m = V^{S_m}_{K,q}$,  
$\psi_m = 2mV_m + a_m \log(1 + | \cdot |^2)$, and $\eta_m =
    \psi_m - 2\log(1 + | \cdot |^2)$.
We note that $V^{S_m*}_{K,q} \leq V^{S}_{K,q} \leq q$ on $K$, so $V_m$
is upper semicontinuous, and therefore plurisubharmonic. 
By Theorem $4.2.6$ in Hörmander \cite{Hormander:convexity} there exists a
solution $u_m$ to the equation 
   $\bar \partial u_m = f \bar \partial \chi$ such that \eqref{eq:3.2}
   holds.   
    By Corollary $5.4$ in Sigurðardóttir 
\cite{Sig:2024} we have that \eqref{eq:3.3} holds.
By Proposition 4.5 in \cite{MagSigSigSno:2023} 
we have $V_m\leq V^S_{K,q}\leq H_S+c$, for some constant $c$,  
so by Theorem 7.2  in \cite{MagSigSigSno:2023}  it follows  
that 
$p_m \in \mathcal{P}^{\widehat S_{\Gamma_m}}_m(\mathbb{C}^n)$.

We turn our attention to finding an upper bound for $\|(f - p_m)e^{-mq}\|_K$.
To this end  we take  $\varepsilon > 0$, 
$\gamma\in]0,R[$, and  $m_0$ such that $V_m > V^S_{K,q} - \varepsilon$,
on the compact set 
$\overline X_R + \overline{B}(0, 1)$, for $m \geq m_0$.
By the continuity of $V^S_{K,q}$ and the compactness of $K$, we can take 
$\delta \in ]0, 1]$ small enough that
    $B(z, \delta) \subset X_{R - \gamma}$, for all $z \in K$, and
    $V^S_{K,q}(w) < q(z)+\varepsilon/2$ for all $w\in B(z,\delta)$ and
$z \in K$.
Since $\chi|_K = 1$, we have that $p_m = f - u_m$, on $K$.
So, by the Cauchy-Schwarz inequality,
\begin{align*}
    |f(z) - p_m(z)|
    &=
    |u_m(z)|
    \leq
    \Omega^{-1}_{2n} \delta^{-2n} \int_{B(z, \delta)} |u_m|\, d\lambda
\\ 
    &\leq
    \Omega^{-1}_{2n} \delta^{-2n}
    \bigg (
        \int_{B(z, \delta)} |u_m|^2 e^{-\psi_m}\, d\lambda
        \int_{B(z, \delta)} e^{\psi_m}\, d\lambda\bigg )^{1/2},
    \quad
    z \in K,
\end{align*}
where $\Omega_{2n}$ denotes the volume of the unit ball in $\mathbb{R}^{2n}$.
Since $0 < a_m < 1$, we have
\begin{align*}
    \int_{B(z, \delta)} e^{\psi_m}\, d\lambda
    &\leq
    \int_{B(z, \delta)} e^{2m V^S_{K,q}} (1 + | \cdot |^2)\, d\lambda
\\
    &\leq
    e^{2mq(z)+m\varepsilon} \int_{X_R} (1 + | \cdot |^2)\, d\lambda,
    \quad
    z \in K,
\end{align*}
and by \eqref{eq:3.2}, we have
\begin{align*}
    \int_{B(z, \delta)} |u_m|^2 e^{-\psi_m}\, d\lambda
    &\leq
    \frac{1}{a_m} \int_{M_\gamma} |f \bar \partial \chi|^2 e^{-\eta_m}\, d\lambda
\\
    &\leq
    \frac{e^{m \varepsilon}}{a_m(R - \gamma)^{2m}}
        \int_{M_\gamma} |f \bar \partial \chi|^2 (1 + | \cdot |^2)^2\, d\lambda,
    \quad
    z \in K,
\end{align*}
where $M_\gamma = \supp \bar \partial \chi$.
The last step follows from the fact that $\chi|_{X_{R - \gamma}} = 1$,
    so $V^S_{K,q} > \log(R - \gamma)$ on $M_\gamma$, and thus
    $V_m > \log(R - \gamma) - \varepsilon$.
Combining these inequalities we have that
\begin{equation*}
    |f(z) - p_m(z)|e^{-mq(z)}
    \leq
    \frac{C_{\varepsilon, \gamma} e^{m \varepsilon}}{a_m^{1/2} (R - \gamma)^m},
    \quad
    m \geq m_0,
\end{equation*}
where $C_{\varepsilon, \gamma}$ is a constant that does not depend on $m$,
    and finally, by \eqref{eq:3.3},
\begin{equation*}
    \varlimsup_{m \rightarrow \infty}
        d^S_{K, m}(f)^{1/m}
    \leq
    \varlimsup_{m \rightarrow \infty}
        \frac{C_{\varepsilon, \gamma}^{1/m} e^\varepsilon}{a_m^{1/(2m)} (R - \gamma)}
    =
    \frac{e^\varepsilon}{R - \gamma}.
\end{equation*}
Since $\varepsilon > 0$ and $\gamma \in ]0, R[$ are arbitrary
the estimate \eqref{eq:3.1} follows.
\end{prooftx}

\section{Distances from $mS$ to the integer lattice}
\label{sec:04}
The presence of $S_m$ in the proof of Theorem \ref{thm:1.2} is to obtain control of the constants
    $a_m$ in \eqref{eq:3.2}.
If $S$ was used instead of $S_m$ we would need to study
\begin{equation}
\label{eq:4.1}
    a
    =
    \varliminf_{m \rightarrow \infty}
        d(mS, \mathbb{N}^n \setminus mS)^{1/m}.
\end{equation}
The distances $d(mS, \mathbb{N}^n \setminus mS)$ come from Theorems 3.6 and 7.2 in \cite{MagSigSigSno:2023},
    which give sufficient conditions for an entire function to belong to certain polynomial classes,
    depending on $S$.
Lemma 5.3 in \cite{Sig:2024} gives the lower bound
\begin{equation*}
    d(S, \mathbb{N}^n \setminus S)
    \geq
    \frac{1}{\sqrt{n} (n - 1)! M^{n - 1}}
\end{equation*}
when $S$ is a polytope with vertices in $\mathbb{N}^n \cap [0, M]^n$.
In this section we will find a lower bound for $d(mS, \mathbb{N}^n \setminus mS)$,
    that depends on $S$ but not $m$,
    when $S$ is polytope with vertices in $\mathbb{Q}^n_+$.

\smallskip
Recall that for every convex set $A \subset \R^n$ 
and $m\in \N^*$ we have that $mA=\sum_{j = 1}^n A = A+\cdots+A$ with $m$ terms in
the right hand side. Denoting the extremal set of $A$ by $\ext\, A$ and
with $x\in (n+1)A$, the
Min\-kowski theorem \cite[Theorem 2.1.9]{Hormander:convexity}
tells us that we can find $a_0,\dots,a_n$ in $\ext\, A$
such that
\begin{equation*}
  x=\lambda_0a_0+\cdots+\lambda_na_n, \quad \lambda_j\geq 0, \
  j=0,\dots, n, \quad 
  \sum_{k=0}^n\lambda_k=n+1. 
\end{equation*}
By renumbering $a_0,\dots,a_n$ we may assume that $\lambda_0\geq 1$,
which implies that we can write $x=t+a_0$, where
$t=(\lambda_0-1)a_0+\lambda_1a_1+\cdots+\lambda_na_n$ and $a_0\in \ext\, A$.  
By induction 
\begin{equation*}
  mA=nA+\sum_{j=1}^{m-n}\ext\, A, \quad m>n,
\end{equation*}
and if we set $T=(1/n)A$, then
\begin{equation}
  \label{eq:5.1}
  mA=mnT=nT+\sum_{j=1}^{mn-n} \ext\, T
=A+\dfrac 1n \sum_{j=1}^{mn-n}\ext\,  A, \quad m>1.
\end{equation}

\begin{proposition}
\label{Prop:4.1}
Let $S=\ch \{v_1,\dots,v_N\}$
be a convex polytope in $\R^n_+$ with
rational vertices $v_j\in \Q^n_+$.
Then
\begin{equation}
\label{eq:5.2}
    d(mS,\N^n\setminus mS)
    \geq
    \frac{1}{nq} d(nqS, \mathbb{N}^n \setminus nqS),
    \quad
    m\in\N^*,
\end{equation}
where $q$ is the common denominator for all the coordinates of $v_1, \dots, v_n$.
\end{proposition}

\begin{proof}
Let $d_m=d(mS,\N^n\setminus mS)$.
Then $v_j\in (1/q)\N^n$ for $j=1,\dots,N$.
Let $s\in mS$ and $u\in \N^n\setminus mS$ such that
$d_m=d(s,u)$.    By \eqref{eq:5.1} we have
$s=t+h$ where $t\in S$ and 
$$
h\in \dfrac 1n\sum_{j=1}^{mn-n}\ext\, S\subset \dfrac 1{nq} \N^n.
$$ 
Since $u\in \N^n\setminus mS\subseteq (1/(nq))\N^n$ we have
$u-h\in (1/(nq))\N^n$.  We have $u-h\not\in S$
for otherwise \eqref{eq:5.1} would imply that 
$u=(u-h)+h\in S$ which does not hold.
Hence 
$$
d_m=d(s,u)=d(s-h,u-h)=d(t,u-h)\geq 
d(S,(1/(nq))\N^n\setminus S),
$$
concluding the proof.
\end{proof}

\medskip
The inclusion of $S_m$ in the proof of Theorem \ref{thm:1.2} would be unnecessary if we could show that $a$
    in \eqref{eq:1.2} was always $1$.
This is not the case.
In fact, we can explicitly construct a lower set $S$ such that $a = 0$.
Let $f \colon [0, 1] \rightarrow [0, 1]$, given by
\begin{equation*}
    f(t)
    =
    1 - e^{-ct^{-b} + c},
    \quad
    t \in [0, 1],
\end{equation*}
where $b > 1$ and $c > 1 + 1/b$, and
\begin{equation*}
    S =
    \{
        x = (x_1, x_2) \in \mathbb{R}^2
        \,;\,
        0 \leq x_2 \leq f(x_1), 0 \leq x_1 \leq 1
    \}.
\end{equation*}
For $t \in ]0, 1[$, we have
\begin{equation*}
    f'(t)
    =
    bct^{-b - 1}(f(t) - 1)
    \quad
    \text{and}
    \quad
    f''(t)
    =
    bct^{-2b - 2}(bc - (b + 1)t^b)(f(t) - 1).
\end{equation*}
Since $f(t) < 1$ for $t \in ]0, 1[$, we have that $f$ is decreasing and concave on $[0, 1]$, so $S$ is convex.
Furthermore, $d(S, (\delta, 1)) \leq d((\delta, f(\delta), (\delta, 1)) = 1 - f(\delta)$, so
\begin{equation*}
    d(mS, \mathbb{N}^n \setminus mS)
    =
    md(S, (1/m)\mathbb{N}^2 \setminus S)
    \leq
    m(1 - f(1/m)).
\end{equation*}
Consequently, we have that
\begin{equation*}
    d(mS, \mathbb{N}^n \setminus mS)^{1/m}
    \leq
    m^{1/m}(1 - f(1/m))^{1/m}
    \leq
    m^{1/m}e^{-cm^{b - 1} + c/m}
    \rightarrow
    0,
\end{equation*}
as $m \rightarrow +\infty$.
So $a$ in \eqref{eq:4.1} is $0$.

{\small 
\bibliographystyle{siam}
\bibliography{rs_bibref}

\smallskip\noindent
Science Institute, University of Iceland,
IS-107 Reykjav\'ik,  ICELAND. 

\smallskip\noindent
bsm@hi.is, ragnar@hi.is, bergur@hi.is.
}

\end{document}